\documentclass{amsart}
\usepackage[english]{babel}
\usepackage[latin1]{inputenc}
\usepackage[dvips,final]{graphics}
\usepackage{amsmath,amsfonts,amssymb,amsthm,amscd,array,
stmaryrd,mathrsfs,mathdots,epigraph}
\usepackage[makeroom]{cancel}
\usepackage{pstricks}
\usepackage[all]{xy}
\usepackage{url}
\usepackage{multirow, blkarray}
\usepackage{booktabs}
\usepackage{caption}
\usepackage{subcaption}

\usepackage{blindtext}
\usepackage{textcomp}
\usepackage[final]{epsfig}
\usepackage{color}
\usepackage{mathabx}

\theoremstyle{plain}
\newtheorem{fac}{Fact}

\theoremstyle{definition}

\newcommand{\R}{\mathbb{R}}

\newcommand{\C}{\mathbb{C}}

\def\om{\omega}
\begin{document}

\title{On Ramanujan's cubic composition formula}

\author{Valentin Ovsienko}
\address{
Valentin Ovsienko,
Centre National de la Recherche Scientifique,
Laboratoire de Math\'ematiques de Reims, UMR9008 CNRS,
Universit\'e de Reims Champagne-Ardenne,
U.F.R. Sciences Exactes et Naturelles,
Moulin de la Housse - BP 1039,
51687 Reims cedex 2,
France}
\email{valentin.ovsienko@univ-reims.fr}

\maketitle

\thispagestyle{empty}

The identity
\begin{equation}
\label{NormEq}
(a_0^2+a_1^2)\,(b_0^2+b_1^2)=
(a_0b_0-a_1b_1)^2+(a_0b_1+a_1b_0)^2,
\end{equation}
is as old as mathematics.
It was known to Diophantus in III century and
appeared again in the early VII century in a book written by Brahmagupta,
it was also used by Fibonacci (who translated Brahmagupta's book from Sanskrit) 
in his ``Book of Squares''.
The identity \eqref{NormEq} is commonly known under the name of 
{\it Brahmagupta-Fibonacci $2$-square identity}.

The identity \eqref{NormEq} holds for elements of an arbitrary
commutative ring, but the case of real coefficients is the most useful.
In this case, a conceptual interpretation of~\eqref{NormEq} uses complex numbers
and is nothing else than the multiplication property of the norm.
Consider the complex numbers:
$$
a=a_0+a_1i,
\qquad
b=b_0+b_1i,
$$
the identity then reads
$$
\left|a\right|^2\left|b\right|^2=
\left|ab\right|^2.
$$
One can argue that, in a sense, the algebra $\C$ of complex numbers was known for more than a thousand
years before its official discovery.

The identity \eqref{NormEq} is the simplest {\it square identity},
i.e., an identity rewriting a product of  sums of squares again
as a sum of squares (with bilinear expressions in the right-hand-side).
The theory of square identities was initiated by Hurwitz~\cite{Hur},
who proved his celebrated $1,2,4,8$ theorem stating that identities with equal number
of summands exist only in dimensions $1,2,4$ and $8$.
Square identities appear, sometimes unexpectedly, 
in many different areas, such as 
algebra and representation theory, geometry and topology, combinatorics, coding theory,
connecting them in a beautiful manner.
For an account on the theory; see the book~\cite{Sha}.
The geometric and topological picture involves
the exceptional spheres $S^0,S^1,S^3$, and~$S^7$.
The theorem of Bott-Milnor~\cite{BM} and Kervaire~\cite{Ker} states that these are the only
parallelizable spheres.
The spheres $S^0,S^1,S^3,S^7$ also appear as fibers of the Hopf fibrations,
(while $S^1,S^2,S^4,S^8$ as their bases).
Some geometric applications are discussed in~\cite{OT}.

The subject historically belongs to number theory, mostly because of the Fermat and Lagrange theorems
about representation of integers as sums of two and four squares, respectively.
The Brahmagupta-Fibonacci $2$-square identity and analogous
Euler's $4$-square identity play crucial role in their proofs.

A natural (yet naive!) question can be asked: what about sums of cubes?
An attempt to find an identity similar to~\eqref{NormEq} representing a product of
sums of cubes fails.
However, a more clever variant of the question was considered and 
answered by no less than Ramanujan.

\newpage 

%%%%%%%%%%%%%%%%%%%%
\section{Ramanujan's cubic identity and its ``demystification''}
%%%%%%%%%%%%%%%%%%%%

The following cubic identity belongs to
Ramanujan; see~\cite{Ber}, p.21.
Instead of a plane sum of cubes, he uses the cubic form
$$
\mathcal{C}(a)=a_0^3+a_1^3+a_2^3-3\,a_0a_1a_2,
$$
which is well-known in algebraic geometry, 
since it defines a singular cubic curve in the projective plane.

\begin{fac}
\label{RamFac}
Every two triplets of real (or complex, or any other commuting) numbers
$a=(a_0,a_1,a_2)$ and $b=(b_0,b_1,b_2)$ satisfy
\begin{equation}
\label{RamEq}
\mathcal{C}(a)\,\mathcal{C}(b)
=
\mathcal{C}(c),
\end{equation}
where
\begin{equation}
\label{BiComEq}
\begin{array}{rcl}
c_0&=&a_0b_0+a_1b_2+a_2b_1,\\[4pt]
c_1&=&a_0b_1+a_1b_0+a_2b_2,\\[4pt]
c_2&=&a_0b_2+a_1b_1+a_2b_0.\\
\end{array}
\end{equation}
\end{fac}

Ramanujan does not reveal a proof, neither he explains how he find
his identity.
Nevertheless, the identity~\eqref{RamEq} and the bilinear composition~\eqref{BiComEq}
are of no mystery.
The simplest way to understand and prove~\eqref{RamEq}
is of course not a direct computation.
Consider $3\times3$ circulant matrices of the form
\begin{equation}
\label{Cirk}
A(a)=
\begin{pmatrix}
a_0&a_1&a_2\\[3pt]
a_2&a_0&a_1\\[3pt]
a_1&a_2&a_0
\end{pmatrix}.
\end{equation}
These matrices form an algebra with the multiplication
$A(a)\,A(b)=A(c)$, where~$c$ is given by~\eqref{BiComEq}.
The cubic form~$\mathcal{C}$ then coincides with the determinant:
$$
\mathcal{C}(a)=\det\left(A(a)\right),
$$
and so Ramanujan's identity immediately follows from the multiplicative property of the determinant.

To compare~\eqref{RamEq} and~\eqref{NormEq}, 
recall that complex numbers can be represented
by $2\times2$~matrices of the form
$$
\begin{pmatrix}
a_0&-a_1\\[3pt]
a_1&a_0
\end{pmatrix},
$$
with $a_0,a_1\in\R$.
The squared norm $|a|^2$ coincides with the determinant.
The main difference between the above real $2\times2$~matrices and $3\times3$~matrices~\eqref{Cirk}
is that the $2\times2$~matrices are always non-degenerate,
unless both~$a_0$ and~$a_1$ vanish.
There is no similar ``full rank property'' in the $3$-dimensional case,
and therefore the algebra of $3\times3$ circulant matrices has zero divisors.

A classification of algebras equipped with invariant forms of arbitrary degree
is known; see~\cite{McC,Sch}.
This result can be understood as analog of Hurwitz's $1,2,4,8$ theorem.
In this sense, cubic identities similar to Ramanujan's
identity may be considered as classified in somewhat non-explicit way.
But this does not detract from elegance of Ramanujan's identity
and nice properties of his cubic form that we now briefly investigate.

%%%%%%%%%%%%%%%%%%%%
\section{Singular but non-degenerate}
%%%%%%%%%%%%%%%%%%%%
We first consider very simple properties of $\mathcal{C}$ that can be easily checked
and leaved to the reader as an exercise.

The cubic form~$\mathcal{C}$ is {\it singular}.
Indeed, it factorizes  
$$
\mathcal{C}(a)=
\left(a_0+a_1+a_2\right)
\left(a_0^2+ a_1^2+ a_2^2-a_0 a_1-a_0a_2-a_1a_2\right).
$$
Note that over~$\C$ the last term further factorizes as as a product
of two linear forms
$$
\left(a_0+\om a_1+\om^2 a_2\right)
\left(a_0+\om^2 a_1+\om a_2\right),
$$
where~$\om=e^{(2i\pi)/3}$ is the cube root of~$1$.
The complex curve $\mathcal{C}(a)=0$ consists in three projective lines;
its real points contain one real projective line $a_0+a_1+a_2=0$
and one point $a_0=a_1=a_2$.

On the other hand, the cubic form~$\mathcal{C}$ is {\it non-degenerate}.
This means that the corresponding symmetric trilinear form, obtained as a polarization
of~$\mathcal{C}$,
$$
\begin{array}{rcl}
(a,b,c) &:=&
\mathcal{C}(a+b+c)\\[4pt]
&&-\mathcal{C}(a+b)-\mathcal{C}(a+c)-\mathcal{C}(b+c)\\[4pt]
&&+\mathcal{C}(a)+\mathcal{C}(b)+\mathcal{C}(c)
\end{array}
$$
satisfies the property of non-degeneressance:
if for a fixed $a$, one has $(a,b,c)=0$ for all~$b$ and~$c$, then $a=0$.

%%%%%%%%%%%%%%%%%%%%
\section{Zero divisors}
%%%%%%%%%%%%%%%%%%%%

As mentioned, the algebra of $3\times3$ circulant matrices~\eqref{Cirk} have divisors of zero.
They belong to the $1$-dimensional diagonal subspace defined by $a_0=a_1=a_2$
and the $2$-dimensional subspace $a_0+a_1+a_2=0$ orthogonal to it.
More precisely, one has the following.

\begin{fac}
\label{DZFac}
Two triplets of real numbers
$a=(a_0,a_1,a_2)$ and $b=(b_0,b_1,b_2)$ satisfy
$$
A(a)\,A(b)
=
0,
$$
if and only if
$$
\left\{
\begin{array}{l}
a_0+a_1+a_2=0 ,\\[4pt]
b_0=b_1=b_2
\end{array}
\right.
\qquad
\hbox{or}
\qquad
\left\{
\begin{array}{l}
a_0=a_1=a_2,\\[4pt]
b_0+b_1+b_2=0.
\end{array}
\right.
$$
\end{fac}

To prove this, notice that $A(a)\,A(b)=0$ implies that both matrices, $A(a)$ and $A(b)$,
are degenerated, in other words, $\mathcal{C}(a)=\mathcal{C}(b)=0$.
Real solutions of this equation were discussed in the previous section.
It remains to check that the product $A(a)\,A(b)$ cannot vanish if the
points of $\mathbb{P}^2$ with coordinates~$a$ and~$b$ belong to the same projective line.

%%%%%%%%%%%%%%%%%%%%
\section{The cubic form $\mathcal{C}$ is its own Hessian}
%%%%%%%%%%%%%%%%%%%%

The cubic form $\mathcal{C}$
has one more remarkable analytic property.

Given a function~$f:\R^3\to\R$ with coordinates on~$\R^3$ still denoted by
$(a_0,a_1,a_2)$, the Hessian of~$f$ is again a function, $H_f:\R\to\R$,
given by the $3\times3$ determinant
$$
H(f)=
\left|
\begin{array}{rcl}
\frac{\partial^2 f}{\partial a_0\partial a_0}&
\frac{\partial^2 f}{\partial a_0\partial a_1}&
\frac{\partial^2 f}{\partial a_0\partial a_2}\\[6pt]
\frac{\partial^2 f}{\partial a_1\partial a_0}&
\frac{\partial^2 f}{\partial a_1\partial a_1}&
\frac{\partial^2 f}{\partial a_1\partial a_2}\\[6pt]
\frac{\partial^2 f}{\partial a_2\partial a_0}&
\frac{\partial^2 f}{\partial a_2\partial a_1}&
\frac{\partial^2 f}{\partial a_2\partial a_2}
\end{array}
\right|.
$$
It expresses geometric properties of~$f$ and can be understood as a sort of curvature.

If $f$ is a cubic form, then its Hessian is again a cubic form.
Calculation of the Hessian of a cubic form can be taken as an iterative process:
$$
f\mapsto H(f)\mapsto H(H(f))\mapsto\cdots
$$
Periodic orbits of this iteration are of interest.

\begin{fac}
\label{DZFac}
The form $\mathcal{C}$ is proportional to its own Hessian:
$$
H(\mathcal{C})=-54\,\mathcal{C}.
$$
\end{fac}

The cubic form $\mathcal{C}$ is characterized by this fact
up to a multiple.

We end this short discussion that composition of cubic forms is a deep subject.
We refer to~\cite{Bar} for very general composition laws of binary and ternary cubic forms.
It seems however, that Ramanujan's identity escapes this theory.
At least it does not appear explicitly among other composition laws considered there.

\end{document}